\documentclass{amsart}
\usepackage{amsxtra,amssymb,multicol}
\title{WHAT IS A PERIOD ?}
\author{Stefan M\"uller-Stach}
\begin{document} 
\rightline{\Huge \bf WHAT IS A PERIOD ?}
\vskip1cm
\rightline{\it \huge Stefan M\"uller-Stach}
\vskip1in

\begin{multicols}{2}

Mathematical constants like $\pi$, $e$ and $\gamma$ arise frequently in
number theory and other areas of mathematics and physics.
Mathematicians have long wondered whether such numbers are
irrational or perhaps even transcendental, that is, not
algebraic.  Because they are solutions of polynomial equations
with rational coefficients, algebraic numbers form a countable
subset of the complex numbers.  Therefore, most complex numbers
are transcendental, although, for any given number, it is usually
difficult to figure out whether it is transcendental.

This essay is about the arithmetic notion of {\it periods}, a countable subalgebra $\mathcal{P}$ 
of the complex numbers defined around 1999 by Maxim Kontsevich and Don Zagier \cite{KZ}. 
Periods contain all algebraic numbers but also many other transcendental numbers important for number theory. 
This notion of periods generalizes in algebraic geometry and yields the 
theory of periods and period domains for algebraic varieties, see \cite{CG}. 

Kontsevich and Zagier define periods as those complex 
numbers whose real and imaginary parts are values of absolutely convergent integrals 
$$
p=\int_\Delta \frac{f(x_1,\ldots,x_n)}{g(x_1,\ldots,x_n)} dx_1 \cdots dx_n. 
$$
Here $f$ and $g$ are polynomials with coefficients in $\mathbb{Q}$, and the integration
domain $\Delta \subset \mathbb{R}^n$ is given by polynomial inequalities with rational coefficients. 

Some initial examples of periods are 
$$
\log(n)=\int_1^n \frac{dx}{x} \text{  and  } \pi = \int_{x^2+y^2\le 1} dx dy.
$$
The representation of a period by an integral is not unique, in the sense that
there are many different integrals representing that period. 

The values $\zeta(s)=\sum_n n^{-s}$ of the Riemann zeta function at positive integers $s \ge 2$, and 
their natural generalizations, the {\it multiple zeta values} 
$$
\zeta(s_1,...,s_k)= \sum_{n_1 > \cdots > n_k > 0} n_1^{-s_1}n_2^{-s_2}\cdots n_k^{-s_k},
$$
for integers $s_i \ge 1$ with $s_1 \ge 2$ are very interesting periods. 
Even for the odd zeta-values $\zeta(3)$, $\zeta(5)$, ... only a few results 
about their irrationality are known. 
By work of Ap\'ery, $\zeta(3)$ is irrational. Like all multiple zeta values, 
it can be represented as an {\it iterated integral}
$$
\zeta(3)=\iiint_{0 < x< y < z < 1}  \frac{dxdydz}{(1-x)yz}.
$$

Using notions from algebraic geometry, one can also define periods in a different form as follows. 
Let $X$ be a smooth algebraic variety over $\mathbb{Q}$ of dimension $d$.  Take
a regular algebraic $d$-form $\omega$ on $X$ and a normal crossing
divisor $D$ in $X$; both $\omega$ and $D$ are also defined over $\mathbb{Q}$.  
Then let $\gamma$ be a singular chain on the underlying topological
manifold $X(\mathbb{C})$ with boundary in $D(\mathbb{C})$. The integral 
$$
p=\int_\gamma \omega
$$
is the period of the quadruple $(X,D,\omega,\gamma)$. 

From an even higher viewpoint, periods are matrix coefficients of the {\it period isomorphism} 
$$
H^*_\text{dR}(X,D) \otimes_\mathbb{Q} \mathbb{C} {\buildrel \cong \over \longrightarrow} H^*_\text{sing}(X,D) \otimes_\mathbb{Q} \mathbb{C}
$$
between algebraic de Rham cohomology and singular cohomology after choosing $\mathbb{Q}$-bases in both groups. 
In this case, $X$ need not be smooth and forms need not be of top degree.
Sophisticated arguments show that all three given definitions of periods agree. 

In this setting, we find $2 \pi i$ as the period of $H^1(X)$ with $X={\mathbb P}^1 \setminus \{0,\infty\}$, $D=\emptyset$, 
$\omega=\frac{dx}{x}$ and $\gamma$ the unit circle.
In a similar way, $\log(n)$ is one of the periods of $H^1(X,D)$, where $X={\mathbb P}^1 \setminus \{0,\infty \}$, $D=\{1,n\}$, $\omega= \frac{dx}{x}$ 
and $\gamma=[1,n]$. 
Certain special $\Gamma$-values occur in the {\it Chowla-Selberg formula} for periods of abelian varieties with complex multiplication. 
The {\it Beilinson conjectures}, which extend Dirchlet's class number formula in algebraic number theory, 
would imply that leading terms in the Taylor series of {\it $L$-functions} of motives are periods in the extended period ring 
$\widehat{\mathcal{P}}=\mathcal{P}[\frac{1}{\pi}]$, where $\pi$ is inverted. A completely different example comes from quantum field theory, 
where periods arise as values of regularized {\it Feynman amplitudes}. Periods of {\it homotopy groups} are another source of examples. 

In addition to the additivity in the integrand and the integration domain, periods inherit from calculus some well-known
relations: a change of variables formula $\int_\gamma f^* \omega=\int_{f_*\gamma} \omega$, and 
Stokes' formula $\int_\gamma \delta \omega=\int_{\partial \gamma} \omega$. Fubini's theorem gives 
$\mathcal{P}$ a multiplication, hence it becomes a $\mathbb{Q}$-algebra. 

At the time Kontsevich and Zagier formulated their idea, not a single explicit non-period number was known. 
In 2008, Masahiko Yoshinaga (arXiv:0805.0349) wrote down a computable non-period, using a variant of Cantor's diagonal argument. Moreover, he showed 
that all periods are elementary computable, i.e., they lie in a certain proper subset of all computable complex numbers, 
so that there are computable non-periods. 

It is still unknown whether $e$ or $1/\pi$ are periods. Presumably they are not. The notion of  
{\it exponential periods} was invented to extend periods to a larger set containing $e$ \cite{KZ}. 

Let us now turn to deeper properties of periods, so that we find out more about the structure of $\mathcal{P}$. 
It turns out that the very abstract viewpoint of {\it mixed motives} provides insights, and 
brings into the game a big symmetry group $G$, the {\it motivic Galois group}. 

Pure and mixed motives were envisioned by Alexander Grothendieck in order to formalize 
properties of algebraic varieties.
In the 1990s, Madhav Nori defined an abelian category $MM(\mathbb{Q})$ of mixed motives over $\mathbb{Q}$. 
In Nori's construction, one starts with a directed graph, where the vertices are pairs of algebraic varieties 
$(X,D)$ defined over $\mathbb{Q}$, and the edges between them are deduced from morphisms of pairs 
$(X,D) \to (X',D')$ (''change of variables'') and chains of inclusions 
$Z \subset D \subset X$ (''Stokes' formula''). The edges thus immediately 
resemble relations among periods, and this is what makes the idea so helpful. 
These arrows are not closed under composition. 
However, if one fixes a representation $T$ with values in vector spaces over $\mathbb{Q}$, e.g. singular or de Rham cohomology, 
then there is a universal diagram category $\mathcal{C}(T)$, and an extension of $T$ as a functor.
After formally inverting the Tate motive $\mathbb{Z}(-1)=(\mathbb{P}^1 \setminus \{0,\infty\},\{1\})$ in $\mathcal{C}(T)$, 
one obtains a $\mathbb{Q}$-linear Tannakian, hence abelian, category $MM(\mathbb{Q})$ without any further assumptions \cite{HMS}. 
The motivic Galois group is the pro-algebraic fundamental group $G=\text{Aut}^\otimes(T)$ of $MM(\mathbb{Q})$ in the Tannaka sense. 
We call $G$ a Galois group, as the viewpoint gives a far-ranging extension of the Galois theory of zero-dimensional varieties.
In $MM(\mathbb{Q})$ cohomology groups of algebraic pairs $(X,D)$ are immediately mixed motives, i.e., 
finite-dimensional $\mathbb{Q}$-representations of $G$, or equivalently comodules over the associated Hopf algebra $A$.  
Both singular and de Rham cohomology provide fiber functors $T_\text{sing}$, $T_\text{dR}$ from $MM(\mathbb{Q})$ to $\mathbb{Q}$-vector spaces. 
The pro-algebraic torsor $\text{Isom}^\otimes(T_\text{dR},T_\text{sing})$ is given by $\text{Spec}(\widehat{\mathcal{P}}_\text{formal})$, 
where $\widehat{\mathcal{P}}_\text{formal}$ is the algebra of formal periods, 
i.e., generated by quadruples $(X,D,\omega,\gamma)$, and subject only to the relations of linearity, change of variables and Stokes \cite{HMS,KZ}.

In this setting, $\widehat{\mathcal{P}}$ is the set of periods of all mixed motives over $\mathbb{Q}$. Multiple zeta values form
the subset of periods of mixed Tate motives over $\mathbb{Z}$. 
The motivic Galois group restricted to mixed Tate motives over $\mathbb{Z}$ gives much control 
over multiple zeta values and implies relations among multiple zeta values $\zeta(s_1,\ldots, s_k)$ 
of a fixed weight $s_1+ \cdots +s_k$. 
The work of Francis Brown on multiple zeta values and the fundamental group of $\mathbb{P}^1 \setminus \{0,1,\infty \}$ 
demonstrates again the value of the philosophy of motives \cite{B}. 
Also in other parts of number theory, for example in the area of rational points, motivic arguments can be applied in finiteness proofs. 

Grothendieck formulated the famous and difficult
{\it period conjecture}, stating that any relation among periods is coming from algebraic geometry, in particular through 
algebraic cycles on products of varieties. In the setting of Nori, this is essentially equivalent to saying that the evaluation map 
$\text{ev}: \widehat{\mathcal{P}}_\text{formal} \to \widehat{\mathcal{P}}$ is injective. This conjecture would have strong 
consequences for the transcendence degree of the space of all periods of a given algebraic variety $X$ via the action of $G$. 

\renewcommand*{\refname}{Further Reading:}

\end{multicols}

\end{document}